\documentclass{amsart}
\usepackage{amssymb}
\usepackage{upref}

\usepackage[all,cmtip]{xy}

\usepackage[lite,initials,shortalphabetic]{amsrefs}

\providecommand{\texorpdfstring}[2]{#1}

\makeatletter
  \renewcommand{\p@enumi}{\thesubsection}
\makeatother

\newenvironment{resumeenumerate}[1]
{\begin{enumerate}
 \setcounter{enumi}{#1}
 \addtocounter{enumi}{-1}
}
{\end{enumerate}
}

\newenvironment{lettered}
{\begin{list}{\thelettercounter)}
 {\usecounter{lettercounter}\def\makelabel##1{\hss\llap{##1}}}
}
{\end{list}
}
\newcounter{lettercounter}
\renewcommand{\thelettercounter}{\alph{lettercounter}}

\makeatletter
\newcommand{\emsection}[1]{%
  \par
  \addpenalty\@secpenalty
  \vskip 6 pt plus 9 pt
  \emph{#1.}\nobreak\enspace\ignorespaces
}
\makeatother

\newcommand{\intro}{%
  \goodbreak
  \vskip 6 pt plus 9 pt
}

\numberwithin{equation}{subsection}

\newcommand{\Period}{\rlap{\enspace .}}

\newcommand{\cat}[1]{\boldsymbol{#1}}

\newcommand{\bs}{\boldsymbol}

\newcommand{\RelCat}{\mathbf{RelCat}}

\newcommand{\Rk}{\mathbf{Rk}}

\newcommand{\Cat}{\mathbf{Cat}}

\newcommand{\SCat}{\mathbf{SCat}}

\newcommand{\DK}{\mathbf{DK}}

\newcommand{\simp}{\mathrm{s}}

\DeclareMathOperator{\Rel}{Rel}

\newcommand{\id}{\mathrm{id}}

\newcommand{\iso}{\approx}

\newcommand{\op}{^{\mathrm{op}}}

\begin{document}

\title[Rezk equivalences are exactly DK-equivalences]
{In the category of relative categories the Rezk equivalences are
  exactly the DK-equivalences}

\author{C. Barwick}
\address{Department of Mathematics, Massachusetts Institute of
  Technology, Cambridge, MA 02139}
\email{clarkbar@gmail.com}

\author{D.M. Kan}
\address{Department of Mathematics, Massachusetts Institute of
  Technology, Cambridge, MA 02139}

\date{\today}

\begin{abstract}
  In a previous paper we lifted Charles Rezk's complete Segal model structure
  on the category of simplicial spaces to a Quillen equivalent one on
  the category of ``relative categories'' and our aim in this
  successor paper is to obtain a more explicit description of the weak
  equivalences in this model structure by showing that \emph{these
    weak equivalences are exactly the DK-equivalences}, i.e.\ those
  maps between relative categories which induce a weak equivalence
  between their simplicial localizations.
\end{abstract}

\maketitle


\section{A formulation of the results}
\label{sec:Rslts}

We start with some preliminaries.
\subsection{Relative categories}
\label{sec:RelCat}

As in \cite{BK1} we denote by $\RelCat$ the category of (small)
relative categories and relative functors between them, where by a
\textbf{relative category} we mean a pair $(\cat C, \cat W)$
consisting of a category $\cat C$ and a subcategory $\cat W \subset
\cat C$ which contains all the objects of $\cat C$ and their identity
maps and of which the maps will be referred to as \textbf{weak
  equivalences} and where by a \textbf{relative functor} between two
such relative categories we mean a weak equivalence preserving
functor.

\subsection{Rezk equivalences}
\label{sec:RzkEq}

In \cite{BK1} we lifted Charles Rezk's complete Segal model structure
on the category $\simp\cat S$ of (small) \emph{simplicial spaces}
(i.e.\ bisimplicial sets) to a Quillen equivalent model structure on
the category $\RelCat$ \eqref{sec:RelCat}.  We will refer to the weak
equivalences in both these model structures as \textbf{Rezk
  equivalences} and denote by both
\begin{displaymath}
  \Rk\subset \simp\cat S
  \qquad\text{and}\qquad
  \Rk\subset\RelCat
\end{displaymath}
the subcategories consisting of these Rezk equivalences.

\subsection{DK-equivalences}
\label{sec:DKeq}

A map in the category $\SCat$ of (small) \emph{simplicial categories}
(i.e.\ categories enriched over simplicial sets) is \cite{Be1} called
a DK-\textbf{equivalence} if it induces \emph{weak equivalences}
between the simplicial sets involved and an \emph{equivalence of
  categories} between their \emph{homotopy categories}, i.e.\ the
categories obtained from them by replacing each simplicial set by the
set of its components.

Furthermore a map in $\RelCat$ will similarly be called a
DK\textbf{equivalence} if its image in $\SCat$ is so under the
\emph{hammock localization functor} \cite{DK2}
\begin{displaymath}
  L^{H}\colon \RelCat\longrightarrow \SCat
\end{displaymath}
(or of course the naturally DK-equivalent functors $\RelCat \to \SCat$
considered in \cite{DK1} and \cite{DHKS}*{35.6}).

We will denote by both
\begin{displaymath}
  \DK\subset \SCat
  \qquad\text{and}\qquad
  \DK\subset \RelCat
\end{displaymath}
the subcategories consisting of these DK-equivalences.

\intro
Our main result then is
\subsection{Theorem}
\label{sec:ThmRzDK}

\emph{A map in $\RelCat$ \eqref{sec:RelCat} is a Rezk equivalence
  \eqref{sec:RzkEq} iff it is a DK-equivalence \eqref{sec:DKeq}.}

\section{An outline of the proof}
\label{sec:Outline}

The proof of theorem~\ref{sec:ThmRzDK} heavily involves the notion of
\subsection{Homotopy equivalences between relative categories}
\label{sec:HomEqRlCt}

A relative functor $f\colon \cat X \to \cat Y$ between two relative
categories \eqref{sec:RelCat} is called a \textbf{homotopy
  equivalence} if there exists a relative functor $g\colon \cat Y \to
\cat X$ (called a \textbf{homotopy inverse} of $f$) such that the
compositions $gf$ and $fg$ are naturally weakly equivalent (i.e.\ can
be connected by a finite zigzag of natural weak equivalences) to the
identity functors of $\cat X$ and $\cat Y$ respectively.
This definition readily implies:
\begin{em}
  \begin{itemize}
  \item [($*$)] If $f\colon \cat X \to \cat Y$ is a homotopy
    equivalence between relative categories which have the two out of
    three property, then a map $x\colon X_{1} \to X_{2} \in \cat X$ is
    a weak equivalence \eqref{sec:RelCat} iff the induced map
    $fx\colon fX_{1} \to fX_{2} \in \cat Y$ is so.
  \end{itemize}
\end{em}

\intro
We also need the following three results.
\subsection{The simplicial nerve functor \texorpdfstring{$N$}{N}}
\label{sec:SimpNv}
In view of \cite{BK1}*{6.1}
\begin{em}
  \begin{itemize}
  \item [($*$)] the simplicial nerve functor $N\colon \RelCat \to
    \simp\cat S$ \eqref{sec:SmpNvFn} is a homotopy equivalence
    \eqref{sec:RzkEq}
    \begin{displaymath}
      N\colon (\RelCat, \Rk)\longrightarrow (\simp\cat S, \Rk) \Period
    \end{displaymath}
  \end{itemize}
\end{em}

\subsection{The relativization functor \texorpdfstring{$\Rel$}{Rel}}
\label{sec:RlFnc}

It was shown in \cite{BK2}*{1.6} that the relativization functor
$\Rel\colon \SCat \to \RelCat$ of the delocalization theorem of
\cite{DK3}*{2.5} is a homotopy inverse of the simplicial localization
functors $(\RelCat, \DK) \to (\SCat, \DK)$ mentioned in
\ref{sec:DKeq}.  This clearly implies that
\begin{em}
  \begin{itemize}
  \item [($*$)] the relativization functor $\Rel\colon \SCat \to
    \RelCat$ \eqref{sec:RlFunc} is a homotopy equivalence
    \begin{displaymath}
      \Rel\colon (\SCat, \DK) \longrightarrow (\RelCat, \DK) \Period
    \end{displaymath}
  \end{itemize}
\end{em}

\subsection{The flipped nerve functor \texorpdfstring{$Z$}{Z}}
\label{sec:FlpNv}

In view of \cite{Be2}*{6.3 and 8.6}
\begin{em}
  \begin{itemize}
  \item [($*$)] the flipped nerve functor $Z\colon \SCat \to \simp\cat
    S$ \eqref{sec:FlipNrv} is a homotopy equivalence
    \begin{displaymath}
      Z\colon (\SCat, \DK) \longrightarrow (\simp\cat S, \Rk) \Period
    \end{displaymath}
  \end{itemize}
\end{em}

\intro
These four results marked ($*$) strongly suggest that theorem
\ref{sec:ThmRzDK} should be true.  To obtain a proof one just has to
show that the functors $N\Rel$ and $Z\colon \SCat \to \simp\cat S$ are
naturally Rezk equivalent.

This will be done in \S \ref{sec:CmplPrf} below.  In fact we will
prove the following somewhat stronger result:
\subsection{Proposition}
\label{prop:FnReEq}

\begin{em}
  The functors
  \begin{displaymath}
    N\Rel \text{ and } Z\colon \SCat\longrightarrow \simp\cat S
  \end{displaymath}
  are naturally Reedy equivalent.
\end{em}

\section{Completion of the proof}
\label{sec:CmplPrf}

Before completing the proof of theorem \ref{sec:ThmRzDK}, i.e.\
proving proposition \ref{prop:FnReEq}, we recall first some of the
notions involved.

\subsection{The simplicial nerve functor \texorpdfstring{$N$}{N}}
\label{sec:SmpNvFn}

This is the functor $N\colon \RelCat \to \simp\cat S$ which sends an
object $X \in \RelCat$ to the bisimplicial set which has as its
$(p,q)$-bisimplices ($p,q \ge 0$) the maps
\begin{displaymath}
  \hat p\times \check q \longrightarrow X \in \RelCat
\end{displaymath}
where $\hat p$ denotes the category $0 \to \cdots \to p$ in which only
the identity maps are weak equivalences and $\check q$ denotes the
category $0 \to \cdots \to q$ in which all maps are weak
equivalences.

\subsection{The relativization functor \texorpdfstring{$\Rel$}{Rel}}
\label{sec:RlFunc}

This is the functor $\Rel\colon \SCat \to \simp\cat S$ which sends an
object $\cat A \in \SCat$ to the pair $(b\cat A, b\bs{\id}) \in
\RelCat$ in which
\begin{enumerate}
\item \label{bcata} $b\cat A$ is the Grothendieck construction on
  $\cat A$ obtained by considering $\cat A$ as a simplicial diagram of
  categories, i.e.\ the category
  \begin{enumerate}
  \item[(i)$'$] which has as objects the pairs $(p,A)$ consisting of
    an integer $p\geq 0$ and an object $A \in \cat A$,
  \item[(i)$''$] which has as maps $(p_{1},A_{1}) \to (p_{2},A_{2})$
    the pairs $(t,a)$ consisting of a simplicial operator $t$ from
    dimension $p_{1}$ to dimension $p_{2}$ and a map $a\colon A_{1}
    \to A_{2} \in \cat A_{p_{2}}$, and
  \item[(i)$'''$] in which the composition is given by the formula
    \begin{displaymath}
      (t',a')(t,a) = \bigl(t't, a'(ta)\bigr)
    \end{displaymath}
  \end{enumerate}
\end{enumerate}
and in which
\begin{resumeenumerate}{2}
\item $\bs{\id}$ denotes the subobject of $\cat A$ consisting of the
  maps of the form $(t,\id)$.
\end{resumeenumerate}

\subsection{The flipped nerve functor \texorpdfstring{$Z$}{Z}}
\label{sec:FlipNrv}

This is the functor $Z\colon \SCat \to \simp\cat S$ which sends an
object $\cat A \in \SCat$ to the simplicial space $Z\cat A$ of which
the space in dimension $k \ge 0$ is the simplicial set $(ZA)_{k}$
which is the disjoint union, taken over all ordered sequences $A_{0},
\ldots, A_{k}$ of objects of $\cat A$, of the products
\begin{displaymath}
  \hom(A_{0},A_{1})\times \cdots\times \hom(A_{k-1}, A_{k}) \Period
\end{displaymath}

\subsection{The opposite \texorpdfstring{$\Gamma\op$}{} of the
  category of simplices functor \texorpdfstring{$\Gamma$}{}}
\label{sec:CtSmpFnc}

This is the functor $\Gamma\op\colon \cat S \to \Cat$ which sends a
simplicial set $X \in \cat S$ to its category of simplices, i.e.\ the
category which has
\begin{enumerate}
  \item as objects the pair $(p,x)$ consisting of an integer $p \ge 0$
    and a $p$-simplex of $X$, and
  \item as maps $(p_{1}, x_{1}) \to (p_{2}, x_{2})$ the simplicial
    operators $t$ from dimension $p_{1}$ to dimension $p_{2}$ such
    that $tx_{1} = x_{2}$.
\end{enumerate}

\intro
We also need
\subsection{Some auxiliary notions}
\label{sec:AuxNtn}

For every object $\cat A \in \SCat$, denote
\begin{itemize}
\item by $Y\cat A$ the simplicial diagram of categories of which the
  category $(Y\cat A)_{k}$ in dimension $k \ge 0$ has as objects the
  sequences of maps in $b\cat A$ \eqref{sec:RlFunc} of the form
  \begin{displaymath}
    \xymatrix{
      {(p_{0},A_{0})} \ar[r]^-{(t_{1},a_{1})}
      & {\cdots} \ar[r]^-{(t_{k},a_{k})}
      & {(p_{k},A_{k})}
    }
  \end{displaymath}
  and as maps the commutative diagrams in $b\cat A$ of the form
  \begin{displaymath}
    \xymatrix{
      {(p_{0},A_{0})} \ar[r]^-{(t_{1},a_{1})} \ar[d]_{(u_{0},\id)}
      & {\cdots} \ar[r]^-{(t_{k},a_{k})}
      & {(p_{k},A_{k})} \ar[d]^{(u_{k},\id)}\\
      {(p'_{0},A_{0})} \ar[r]^-{(t'_{1},a'_{1})}
      & {\cdots} \ar[r]^-{(t'_{k},a'_{k})}
      & {(p'_{k},A_{k})}
    }
  \end{displaymath}
\end{itemize}
and
\begin{itemize}
\item by $\overline{Y\cat A} \subset Y\cat A$ the subobject of which
  the category $(\overline{YA})_{k}$ in dimension $k$ is the
  subcategory of $(YA)_{k}$ consisting of the above maps for which
  \emph{the $t_{i}$'s and the $t'_{i}$'s are identities and hence all
    $p_{i}$'s are the same, all $p'_{i}$'s are the same and all
    $u_{i}$'s are the same}.
\end{itemize}

Then there is a \emph{strong deformation retraction} of $Y\cat A$ onto
$\overline{Y\cat A}$ which to each object of $Y\cat A$ as above
assigns the map
\begin{displaymath}
  \xymatrix@C=5em{
    {(p_{0},A_{0})} \ar[r]^-{(t_{1},a_{1})}
    \ar[d]_{(t_{k}\cdots t_{1},\id)}
    & {\cdots} \ar[r]^-{(t_{k},a_{k})}
    & {(p_{k},A_{k})} \ar[d]^{(\id,\id)}\\
    {(p_{k},A_{0})} \ar[r]^-{(\id,t_{k}\cdots t_{1}a_{1})}
    & {\cdots} \ar[r]^-{(\id,a_{k})}
    & {(p_{k},A_{k})}
  }
\end{displaymath}
the existence of which implies that
\begin{em}
  \begin{enumerate}
  \item the inclusion $\overline{Y\cat A} \subset Y\cat A$ is a
    dimensionwise weak equivalence of categories.
  \end{enumerate}
\end{em}

One also readily verifies that there is a canonical 1-1 correspondence
between the objects of $(\overline{Y\cat A})_{k}$ and the simplices of
$(Z\cat A)_{k}$ \eqref{sec:FlipNrv} and that in effect
\begin{em}
  \begin{resumeenumerate}{2}
  \item there is a canonical isomorphism $\overline{Y\cat A} \iso
    \Gamma\op Z\cat A$ \eqref{sec:CtSmpFnc}.
  \end{resumeenumerate}
\end{em}

\intro
Now we are ready for the
\subsection{Completion of the proof}
\label{sec:Cmpltn}

Let $n\colon \Cat \to \cat S$ denote the \emph{classical nerve
  functor}.

Then clearly $N\Rel\cat A = nY\cat A$ \eqref{sec:AuxNtn} and if one
defines $\overline{N\Rel\cat A} \subset N\Rel\cat A$ by
$\overline{N\Rel\cat A} \subset n\overline{Y\cat A}$, then it follows
from \ref{sec:AuxNtn}(i) that
\begin{em}
  \begin{enumerate}
  \item the inclusion $\overline{N\Rel\cat A} \to N\Rel \cat A \in
    \simp\cat S$ is a Reedy equivalence.
  \end{enumerate}
\end{em}
Moreover it follows from \ref{sec:AuxNtn}(ii) that
\begin{em}
  \begin{resumeenumerate}{2}
  \item there is a canonical isomorphism $\overline{N\Rel\cat A} \iso
    n\Gamma\op Z\cat A$
  \end{resumeenumerate}
\end{em}
and to complete the proof of proposition \ref{prop:FnReEq} and hence
of theorem \ref{sec:ThmRzDK} it thus suffices, in view of the fact
that clearly
\begin{em}
  \begin{resumeenumerate}{3}
  \item the functors $n\Gamma\op$ and $n\Gamma\colon \cat S \to \cat
    S$ \eqref{sec:CtSmpFnc} are naturally weakly equivalent,
  \end{resumeenumerate}
\end{em}
to show that
\begin{em}
  \begin{resumeenumerate}{4}
  \item there exists a natural Reedy equivalence $n\Gamma Z\cat A \to
    Z\cat A \in \simp\cat S$.
  \end{resumeenumerate}
\end{em}
But this follows immediately from the observation of Dana Latch
\cite{L} (see also \cite{H}*{18.9.3}) that there exists a natural weak
equivalence
\begin{displaymath}
  n\Gamma \longrightarrow 1 \Period
\end{displaymath}


\end{document}